\documentclass[twoside,leqno]{article}
\usepackage{amsmath,jamsl2e}
\usepackage{latexsym}
\input{epsf.sty}

\newcommand\ack{\section*{Acknowledgement}}

\newcommand\ds{\displaystyle}

\newenvironment{figMacPc}[4]				
{\begin{figure}[t]
	\epsfysize=#2\centerline{\epsfbox{#1}}		
	\caption{#3}					
	\label{#4}							
}
{\end{figure}}
\newcommand{\fine}{\quad\rule{2mm}{2mm}}                    
\begin{document}
\title
{TOWARDS THE MODELING OF NEURONAL FIRING BY GAUSSIAN PROCESSES}
%
%
%
\author
{\sc E. Di Nardo, A.G. Nobile, E. Pirozzi and L.M. Ricciardi}

\date{}


\maketitle

\par\noindent
\begin{abstract}
This paper focuses on the outline of some computational 
methods for the approximate solution of the  integral equations
for the neuronal firing probability density and an algorithm for the generation of
sample-paths in order to construct histograms
estimating the firing densities. Our results originate from
the study of non-Markov stationary Gaussian neuronal models
with the aim to determine the neuron's firing probability density function.
A parallel algorithm has been implemented in order to simulate  large
numbers of sample paths of Gaussian processes characterized by damped
oscillatory covariances in the presence of time dependent
boundaries. The analysis based on the simulation procedure 
provides an alternative research tool when closed-form results or analytic
evaluation of the neuronal firing densities are not available.
\end{abstract}
%

%

%
\section{Introduction}
\hfill\\
This contribution deals with the implementation of procedures and
methods, worked out in our group during the last few years, in
order to provide algorithmic solutions to the problem of determining
the first passage time (FPT) probability densities (pdf) and its
relevant statistics for continuous state-space and continuous
parameter Gaussian processes describing the stochastic modeling
of a single neuron's activity.
\par
In most modeling approaches, it is customary to assume that a neuron 
is subject
to input pulses occurring randomly in time, (see, for instance, \cite{Ricciardi02}
 and references therein). 
As a consequence of the received stimulations, it reacts by 
producing a response that consists of a spike train. The reproduction of the 
statistical features of such spike trains has been the goal of many researches 
who have focused the attention on the analysis of the interspike intervals.
Indeed, the importance of the interspike
intervals is due to the generally accepted hypothesis that
the information transferred within the nervous system is
usually encoded by the timing of occurrence of neuronal spikes.
\par
To describe the dynamics of the neuronal firing we consider 
a stochastic process $X(t)$ representing the change in 
the neuron membrane potential between two consecutive spikes 
(cf., for instance, \cite{Ricciardi77}). In this context, 
the threshold voltage is viewed as a deterministic function $S(t)$ and  
the instant when the membrane potential crosses $S(t)$ as a FPT random variable.
\par
The modeling of a single neuron's activity by means of a
stochastic process has been the object of numerous investigations
during the last four decades. A milestone contribution in this direction is
the much celebrated paper by Gerstein and Mandelbrot
\cite{Gerstein64} in which a random walk and its continuous
diffusion limit (the Wiener process) was proposed with the aim of describing
 a possible, highly schematized,
spike generation mechanism. However, despite the excellent fitting of a variety of data,
this model has been the target of severe criticism on the base of its extreme idealization
in contrast with some electrophysiological evidence: for
example, this model does not take into account the 
spontaneous exponential decay of the neuron membrane potential.
 An improved model is the
so called Ornstein-Uhlenbeck (OU) model, that embodies the
presence of such exponential decay. 
However, the OU model does not allow
to obtain any closed form expression for the firing pdf, except for
some very particular cases of no interest within the neuronal
modeling context. Rather cumbersome computations are thus required to
obtain evaluations of the statistics of the firing time.
Successively, alternative neuronal models have been proposed,
that include more physiologically features.   The literature on this subject is too vast
to be recalled here. We limit ourselves to mentioning that a review of most significant
 neuronal models can be
found in \cite{Ricciardi99}, \cite{Ricciardi02}  and in the references
therein. In particular, in \cite{Ricciardi99} it is presented 
an outline of appropriate mathematical
techniques by which to approach the FPT problem in the neuronal context.
\par
We shall now formally define the firing pdf for a model based on a stochastic
process $X(t)$ with continuous sample paths. First, assume 
$P[X(t_0)=x_0]=1$, with $x_0 < S(t_0),$ i.e. we view the
sample paths of $X(t)$ as originating at a preassigned state $x_0$ at initial time $t_0$.
Then,
$$
T_{x_0} = \inf_{t\geq t_0} \bigl\{t : X(t) > S(t)\bigr\},\qquad x_0<S(t_0)
$$
is the FPT of $X(t)$ through $S(t)$, and
$$
g[S(t),t|x_0,t_0] = \frac{\partial}{\partial t} P(T_{x_0} < t)
$$
is its pdf.
\par
Henceforth, the FPT pdf $g[S(t),t|x_0,t_0]$ will be identified with the firing pdf of a neuron
whose membrane potential is modeled by $X(t)$ and whose firing 
threshold is $S(t)$.
\par

Throughout this paper, we shall focus our attention on neuronal
models rooted on diffusion and Gaussian processes, partially motivated by 
the generally accepted hypothesis
that in numerous instances the neuronal firing
is caused by the superposition of a very large number
of synaptic input pulses which is suggestive of the generation of Gaussian distributions
 by virtue of some sort of central limit theorems. 
\par
It must be explicitly pointed out that models based on diffusion processes
are characterized by the \lq\lq lack of memory\rq\rq as a consequence of 
the underlying Markov property. 
However, in the realistic presence of correlated
input stimulations, the Markov assumption breaks down and one faces the problem of
considering more general stochastic models, for which the literature on FPT
problem appears scarce and fragmentary. 
Simulation procedures then provide possible alternative investigation tools
 especially if they can
be implemented on parallel computers, (see \cite{DiNardo01}).
The goal of  a typical simulation procedure is to sample $N$ values of
the FPT  by a suitable construction of $N$ time-discrete sample
paths of the process and then to record the instants when such
sample paths first cross the boundary. 
In such a way, one is led to obtain estimates of the firing pdf and of its statistics,
that may be implemented for data fitting purposes.
\par
The aim of this paper is to outline numerical and
theoretical methods to characterize the FPT
pdf for Gaussian processes. Attention will be focused on Markov models
in Section \ref{Section2}, and on non-Markov models in Section \ref{Section3}.
Finally, Section \ref{Section4} will be devoted to the description
 of some computational results.
%
\section{Gauss-Markov processes}\label{Section2}
\hfill \\
We start briefly reviewing some essential properties of Gauss-Markov processes.
Let $\{X(t), t\in I\}$, where $I$ is a continuous parameter
set, be a real continuous Gauss-Markov process with the following properties
(cf. \cite{Mehr65}):
\begin{description}
\item[{\it (i)}]  $m(t):= E[X(t)]$ is continuous in $I$;
\item[{\it (ii)}] the covariance $c(s,t):= E\bigl\{[X(s)-m(s)]\,[X(t)-m(t)]\bigr\}$
is continuous in $I\times I$;
\item[{\it (iii)}] $X(t)$ is non-singular, except possibly at the end
points of $I$ where it could be equal to $m(t)$ with probability one.
\end{description}
A Gaussian process is Markov if and only if its covariance
satisfies
\begin{equation}
c(s,u) = {c(s,t)\,c(t,u)\over c(t,t)} \quad \forall s,t,u \in I, s\leq t\leq
u.
\label{eq:(2.1)}
\end{equation}
It is well known  \cite{Mehr65}, that well-behaved solutions of (\ref{eq:(2.1)}) are of the form
\begin{equation}
c(s,t) = h_1(s)\,h_2(t), \qquad s\leq t,
\label{eq:(2.2)}
\end{equation}
where
\begin{equation}
r(t):={h_1(t)\over h_2(t)}
\label{eq:(2.3)}
\end{equation}
is a monotonically increasing function by virtue  of the Cauchy-Schwarz
inequality, and $h_1(t)\,h_2(t)>0$ because of the assumed nonsingularity
of the process on $I$. The conditional pdf $f(x,t|x_0,t_0)$ of
$X(t)$  is  a normal
density characterized respectively by conditional mean and variance
\begin{eqnarray*}
M(t|t_0)& = & m(t)+{h_2(t)\over h_2(t_0)}\;[x_0-m(t_0)]\\
V(t|t_0)& = & h_2(t)\,\biggl[h_1(t)-{h_2(t)\over h_2(t_0)}\;h_1(t_0)\biggr],
\end{eqnarray*}
with $t,t_0\in I,\;t_0<t$. It satisfies the Fokker-Planck equation
and the associated initial condition
\begin{eqnarray*}
& &{\partial f(x,t|y,\tau) \over \partial t}=-\;{\partial \over \partial x}\;
[A_1(x,t)\,f(x,t|y,\tau)]+{1\over 2}\;{\partial^2 \over \partial x^2}\;
[A_2(t)\,f(x,t|y,\tau)],\\
\\
& &\lim_{\tau \uparrow t}\,f(x,t|y,\tau)=\delta(x-y),
\end{eqnarray*}
with $A_1(x,t)$ and $A_2(t)$ given by
$$
A_1(x,t)= m^{\prime}(t) + [x-m(t)]\;{h_2^{\prime}(t)\over h_2(t)}\,,\qquad
A_2(t) = h_2^2(t)\;r^{\prime}(t),
$$
the prime denoting derivative with respect to the argument.
\par
The class of the Gauss-Markov processes $\{X(t), t\in
[0,+\infty)\}$, such that  $f(x,t|y,\tau)\equiv f(x,t-\tau|y),$ is
characterized by means and covariances of the following two forms:
\begin{description}
\item{$\;$} $m(t)=\beta_1 t + c,\quad   c(s,t) = \sigma^2 s + c_1$\\
$\hspace*{4cm}(0 \leq s \leq t < +\infty, \; \beta_1,c\in {\bf R},
\;  c_1\geq 0,\;\sigma\neq 0 )$
\item{or}
\item{$\;$} $m(t)= - {\ds{\beta_1\over\beta_2}} + c\,e^{\beta_2 t},\quad
 c(s,t) = c_1\,e^{\beta_2 t} \left[ c_2\,e^{\beta_2 s} - {\ds{\sigma^2
\over 2 c_1  \beta_2}}\,e^{ - \beta_2 s}\right]$\\
$\hspace*{0.0cm}\left(0 \! \leq \! s \! \leq \! t < \! +\infty, \; \beta_1,c,c_2\in
     {\bf R}, \;\sigma\neq 0,\; c_1\neq 0,\beta_2 \neq 0, \;  c_1c_2
- {\ds{\sigma^2\over 2 \beta_2}}\geq 0\right)$.
\end{description}
The first type includes the Wiener process, used in  \cite{Gerstein64} to describe
the neuronal firing, while  the second type
includes the Ornstein--Uhlenbeck process that has often been invoked as
 a model for neuronal activity (see, for instance, \cite{Ricciardi02}).
\par
Any Gaussian process with covariance as in (\ref{eq:(2.2)})
can be represented in terms of the standard Wiener process $\{W(t),t\geq 0\}$ as
\begin{equation}
X(t)=m(t)+h_2(t)\;W\bigl[r(t)\bigr],
\label{eq:(2.4)}
\end{equation}
and is therefore Markov. This last equation suggests the
way to construct the FPT pdf of a Gauss-Markov process
$X(t)$ in terms of the FPT pdf of the standard Wiener process $W(t).$
As an example, from (\ref{eq:(2.4)}) for the conditioned FPT pdf one has:
\begin{equation}
g[S(t),t|x_0,t_0] ={dr(t)\over
dt}\;g_W\bigl\{S^*[r(t)],r(t)|x_0^*,r(t_0)\bigr\},
\label{eq:(2.12)}
\end{equation}
where $r(t)$ is defined in (\ref{eq:(2.3)}) and
$g_W[S^*(\vartheta),\vartheta|x_0^*,\vartheta_0]$
is the FPT pdf of $W(\vartheta)$ from $x_0^*$ at time $\vartheta_0$ to
the continuous boundary $S^*(\vartheta)$, with
$$x_0^*={x_0-m[r^{-1}(\vartheta_0)]
\over h_2[r^{-1}(\vartheta_0)]},\qquad
S^*(\vartheta)={S[r^{-1}(\vartheta)]-m[r^{-1}(\vartheta)]
\over h_2[r^{-1}(\vartheta)]}\,.$$
Results on the FPT pdf for the standard Wiener process can thus in principle 
be used via
(\ref{eq:(2.12)}) to obtain the FPT pdf of any continuous Gauss-Markov
process.
For instance, if $S^*(\vartheta)$ is linear in $\vartheta$,
$g_W[S^*(\vartheta),\vartheta|x_0^*,\vartheta_0]$ is known and
$g[S(t),t|x_0,t_0]$ can be obtained via (\ref{eq:(2.12)}). Instead, if
$g_W[S^*(\vartheta),\vartheta|x_0^*,\vartheta_0]$ is not known, a
numerical algorithm, or a simulation procedure, should be used
for the standard Wiener process and, after that,
$g[S(t),t|x_0,t_0]$ can be obtained via the indicated transformation. 
\par
The above procedure often exhibits the serious drawback of ensuing unacceptable
time dilations (see \cite{DiNardo01bis}). For instance,  exponentially
large times are involved when transforming the Ornstein-Uhlenbeck
process to the Wiener process, which makes such a method hardly viable.
Hence, it is desirable to dispose of a direct and efficient
computational method to obtain evaluation of the firing pdf. Along such a direction,
in
\cite{DiNardo01bis} it has been proved that the conditioned
FPT density of a Gauss-Markov process can be obtained by solving the 
non-singular Volterra second kind integral equation
\begin{eqnarray}
&&\hspace*{-0.5cm} g[S(t),t|x_0,t_0] \!=\! - 2 \Psi[S(t),t|x_0,t_0] + 2 \int_{t_0}^t\!
g[S(\tau),\tau|x_0,t_0]\,
\Psi[S(t),t|S(\tau),\tau]\;d \tau\nonumber\\
&&\hspace*{10cm}\bigl(x_0<S(t_0)\bigr)
\label{eq:(3.1)}
\end{eqnarray}
with $S(t),m(t),h_1(t), h_2(t) \in C^1(I)$ and
\begin{eqnarray}
&& \Psi[S(t),t\,|\,y,\tau]  =  \biggl\{{S^{\prime}(t) - m^{\prime}(t)\over 2}
  \; -{S(t)-m(t)\over 2}\;{h^{\prime}_1(t)h_2(\tau)
 -h^{\prime}_2(t)h_1(\tau)\over h_1(t)h_2(\tau)-h_2(t)h_1(\tau)} \nonumber \\
&&\hspace*{3.3cm}  - \; {y-m(\tau)\over 2}\;{h^{\prime}_2(t)h_1(t)-h_2(t)
 h^{\prime}_1(t)\over h_1(t)h_2(\tau)-h_2(t)h_1(\tau)}\biggr\}
  \;  f[S(t),t\,|\,y,\tau] \label{psi}
\end{eqnarray}
where $f[x,t|y,\tau]$ is the transition pdf of $X(t).$
Closed form solutions of (\ref{eq:(3.1)}) 
are available in \cite{DiNardo01bis} for different families of
boundaries.
\par
By making use of this result, in \cite{DiNardo01bis} an efficient
numerical procedure based on a repeated Simpson's rule has been proposed
to evaluate FPT densities of Gauss-Markov processes, that can be implemented to obtain reliable
evaluations of firing densities for neuronal models based on Gauss-Markov processes.
%
\section{Gauss non-Markov processes}\label{Section3}
\hfill \\
The methods proposed in the previous Section rest on the strong Markov
assumption on the stochastic process modeling the neuron's membrane potential, 
which grants the possibility of making use of the mentioned
analytic and computational methods for FPT pdf evaluations. 
This is not the case when the stochastic process used to model the neuron's
firing mechanism is non-Markov. Here we shall focus our attention on Gauss non-Markov processes.
 However, thus doing we face the lack of
effective analytical methods for obtaining manageable closed-form
expressions for the FPT pdf, although some
preliminary analytical results have been obtained by Ricciardi and
Sato in \cite{Ricciardi83} for a class of stationary Gaussian processes.
\par
Indeed, if $X(t)$ is one-dimensional non-singular stationary Gaussian process
 mean square differentiable, a series expansion for the FPT pdf was derived
(see, \cite{Ricciardi86}). 
In the following, for convenience, we shall take $t_0=0$ as initial time and, without 
loss of generality, assume that $E[X(t)]=0$ and $P[X(t_0)=x_0]=1$, with $x_0$ an
arbitrarily specified initial state. Furthermore, the covariance function 
$E[X(t)X(\tau )]: =\gamma (t-\tau )$ will be assumed to be  
such that $\gamma (0)=1,\dot \gamma(0)=0$ and $\ddot \gamma (0)<0$
(this last assumptions being equivalent to the mean square differentiable
property). 
The FPT pdf of $X(t)$ through $S(t)$ is then given by
the following expression
\begin{equation}
g[S(t),t|x_0]=W_1(t|x_0)+\sum_{i=1}^\infty
(-1)^i\int_0^t \!\! dt_1 \!\! \int_{t_1}^t \!\! dt_2
\cdots \!\! \int_{t_{i-1}}^t \!\! \!\! dt_i
W_{i+1}(t_1,\ldots ,t_i,t|x_0),  \label{(0)}
\end{equation}
with
\begin{eqnarray}
&&\hspace*{-0.5cm}W_n(t_1, \ldots,  t_n|x_0)  \nonumber\\
&&\hspace*{0.2cm} = 
\int_{\dot S(t_1)}^\infty  dz_1  \cdots
 \int_{\dot S(t_n)}^\infty  dz_n 
\prod_{i=1}^n\,  [z_i-\dot S(t_i)]\,  p_{2n}[  S(t_1),\ldots  ,  S(t_n);
z_1,  \ldots  ,  z_n|x_0],\nonumber
\end{eqnarray}
where $p_{2n}(x_1,\ldots ,x_n;z_1,\! \ldots ,z_n|x_0)$ is the joint pdf of
$\{X(t_1),\ldots \!,X(t_n),Z(t_1)=\dot X(t_1),$ $\ldots, $ $Z(t_n) \! =\dot X(t_n) \}$
conditional upon $X(0)=x_0$. Due to the great complexity of the involved multiple integrals,
 expression (\ref{(0)}) does not appear to be manageable for practical uses, even though it
has recently been shown that it allows to obtain some interesting asymptotic results 
\cite{DiNardo02}.
Since (\ref{(0)}) is a Leibnitz series for each fixed $t>0,$ estimates of the FPT
pdf can in principle be obtained  since its partial sum of order $n$ provides
a lower or an upper bound to $g$ depending on whether $n$ is even or odd.
However, also the evaluation of such partial sums is extremely
cumbersome.
\par
In conclusion, at the present time for this class of Gaussian processes, no
effective analytical methods, nor viable numerical algorithms are available
to evaluate the FPT pdf. A simulation procedure seems to be the only residual way
of approach.
\par
To this aim, we have restored and updated an
algorithm due to Franklin \cite{Franklin65} in order to construct
sample paths of a stationary Gaussian process with spectral
density of a rational type and deterministic starting point.
The idea is the following. Let us consider the linear filter
\begin{equation}
X(t) = \int_{0}^{\infty} h(s)\,W(t-s)\; ds  \label{(3)}
\end{equation}
where $h(t)$ is the impulse response function and $W(t)$ is the input
signal. By Fourier transformation, (\ref{(3)}) yields
\begin{equation}
\Gamma_{X}(\omega) = |H(\omega)|^2\, \Gamma_{W}(\omega) \label{(4)}
\end{equation}
where $\Gamma_{W}(\omega)$ and $\Gamma_{X}(\omega)$ are
the spectral densities of input $W(t)$ and output $X(t)$, respectively, and where
$H(\omega)$ denotes the Fourier transform of $h(t).$
Equation (\ref{(4)}) is suggestive of a method to construct a
Gaussian process $X(t)$ having a preassigned spectral density
$\Gamma_{X}(\omega) \equiv \Gamma(\omega)$.
It is indeed sufficient to consider a white noise $\Lambda(t),$ having
spectral density $\Gamma_{W}(\omega) \equiv 1,$ as the input signal and
then select $h(t)$ in such a way that $|H(\omega)|^2 =\Gamma(\omega)$.
If the spectral density of $X(t)$ is of rational type, namely if
\begin{equation}
\Gamma(\omega) := \int_{-\infty}^{\infty} \gamma(t)\,e^{- i\,\omega\,t}\;dt =
\left| \frac{P(i\,\omega)}{Q(i\,\omega)} \right|^2 \enspace  
\label{(5)}
\end{equation}
where $P$ and $Q$ are polynomials with $deg(P)<deg(Q),$ 
setting $H(\omega)= P(i\, \omega)/Q(i\, \omega)$,
from (\ref{(4)}) it follows
$$X(t)=\frac{P(D)}{Q(D)}\; \Lambda(t)$$
where $D=d/dt.$ To calculate the output signal $X(t)$ it is thus
necessary to solve first the differential equation $Q(D)\,\phi(t)=\Lambda(t)$
to obtain the stationary solution $\phi(t)$, and then evaluate $X(t)=
P(D)\,\phi(t).$
The simulation procedure is designed to construct sample paths of the
process $X(t)$ at the instants $t=0, \Delta t, 2 \Delta t, \ldots$
where $\Delta t$ is a positive constant time increment. The underlying
idea can be applied to any Gaussian process having spectral densities
of a rational type and, since the sample paths of the simulated process
are generated independently of one another, this simulation procedure
is particularly suited for implementation on supercomputers.
\par
Extensive computations have been performed on parallel computers to explore the
different shapes of the FPT pdf as induced by the oscillatory
behaviors of covariances and thresholds (cf., for instance, \cite{DiNardo00} and
\cite{DiNardo00bis}). 
%
\begin{figMacPc}
			{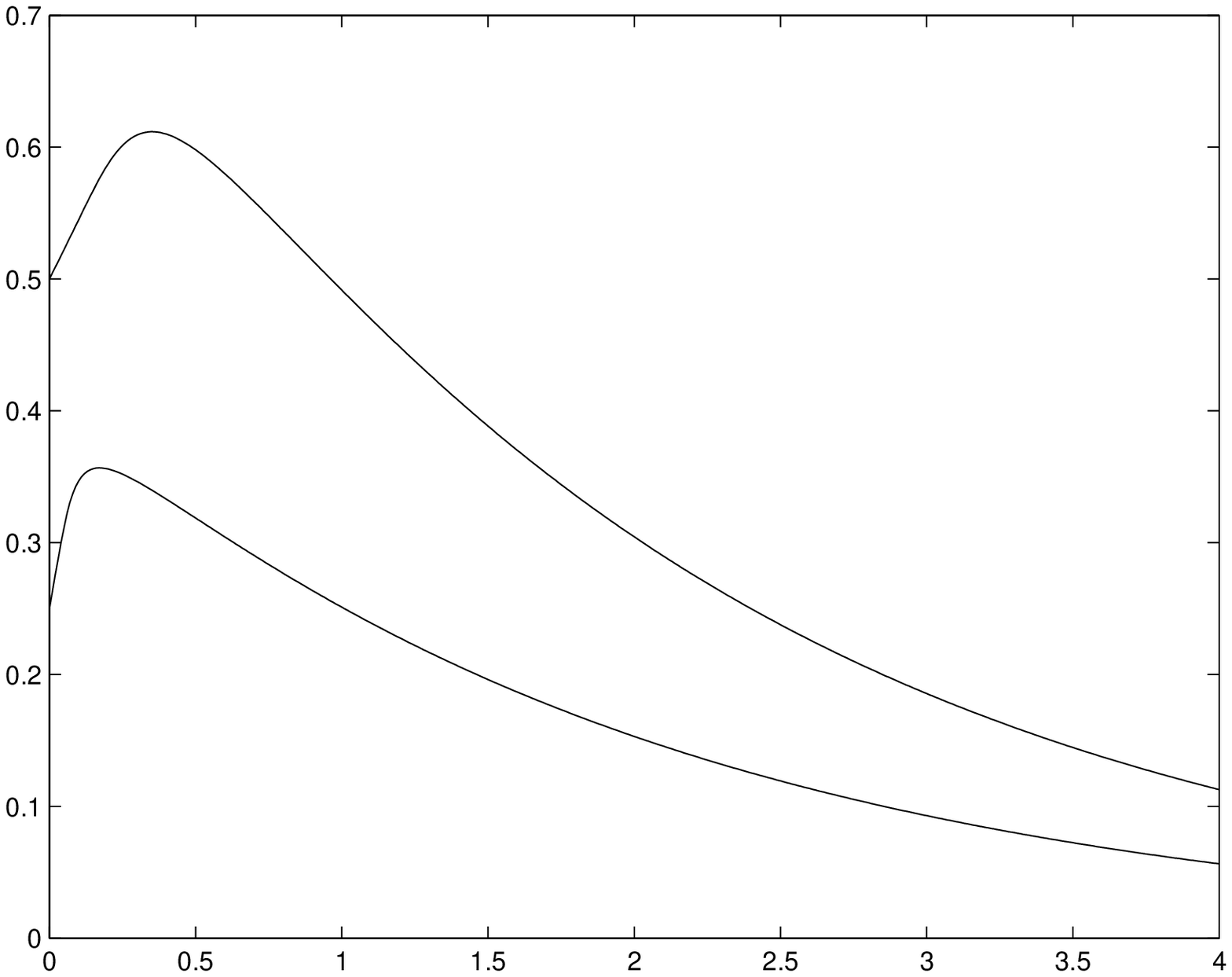}
			{8 cm}
			{{\small Plot of the boundary $S(t)$  given in (\ref{soglia})
for $\beta=0.5$ and $d=0.25, 0.50$ (bottom to top)}}		
	            {fig:fig1}
\end{figMacPc}
%
%
\begin{figMacPc}
			{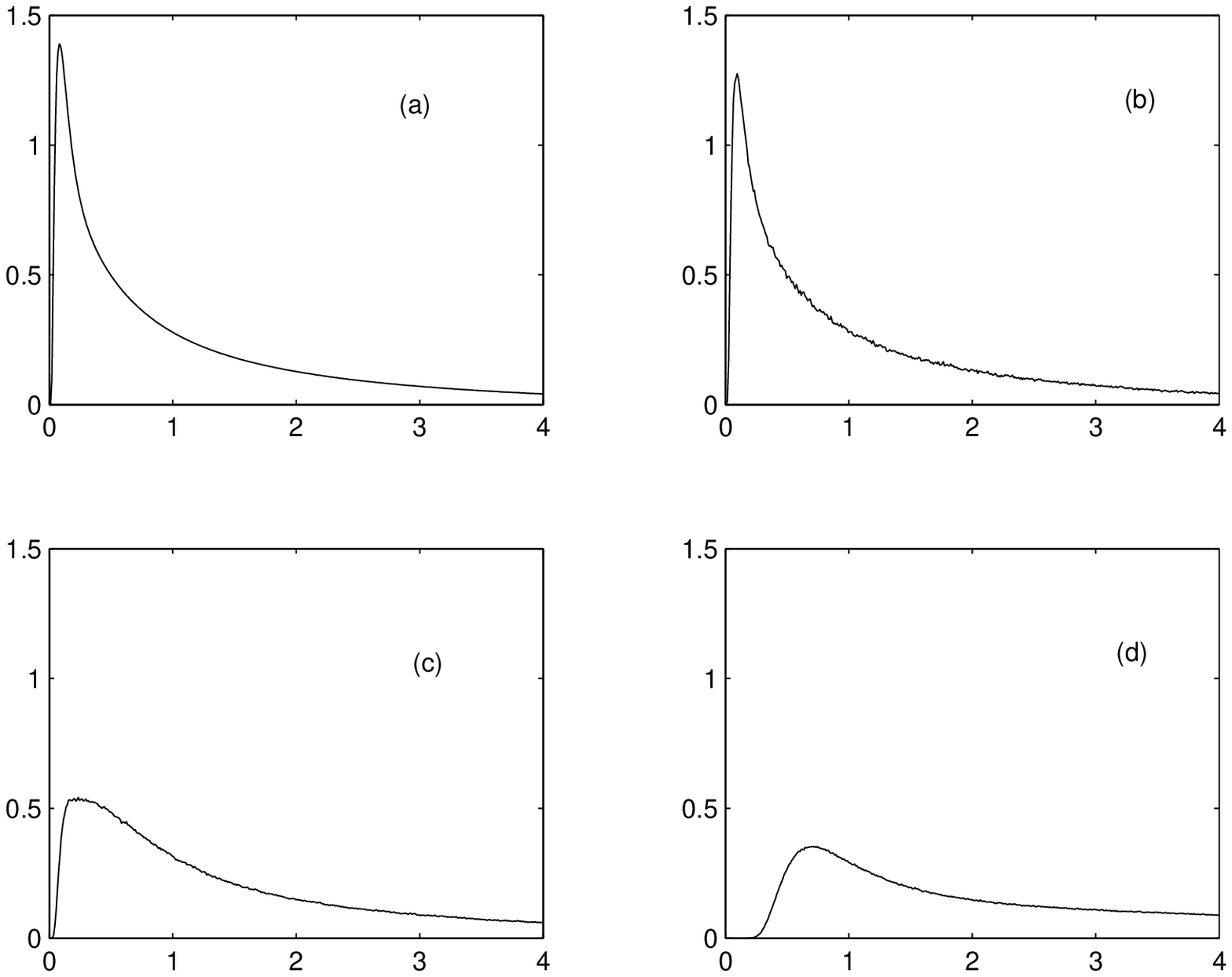}
			{10 cm}
			{{\small Plots refer to  FPT pdf $g(t)$ through the boundary (\ref{soglia})                                                             with $\beta=0.5$ and $d=0.25$ for a zero-mean Gaussian process characterized by correlation function (\ref{(funcorr)}). In Figure~2(a) $g(t)$ given by (\ref{(FPTclosed)}) has been plotted. The estimated FPT pdf $\tilde{g}(t)$ with $\alpha=10^{-10}$ is shown in Figure~2(b), with $\alpha=0.25$ in Figure~2(c) and with $\alpha=0.5$ in Figure~2(d).}}		
	            {fig:fig2}
\end{figMacPc}
%
\begin{figMacPc}
			{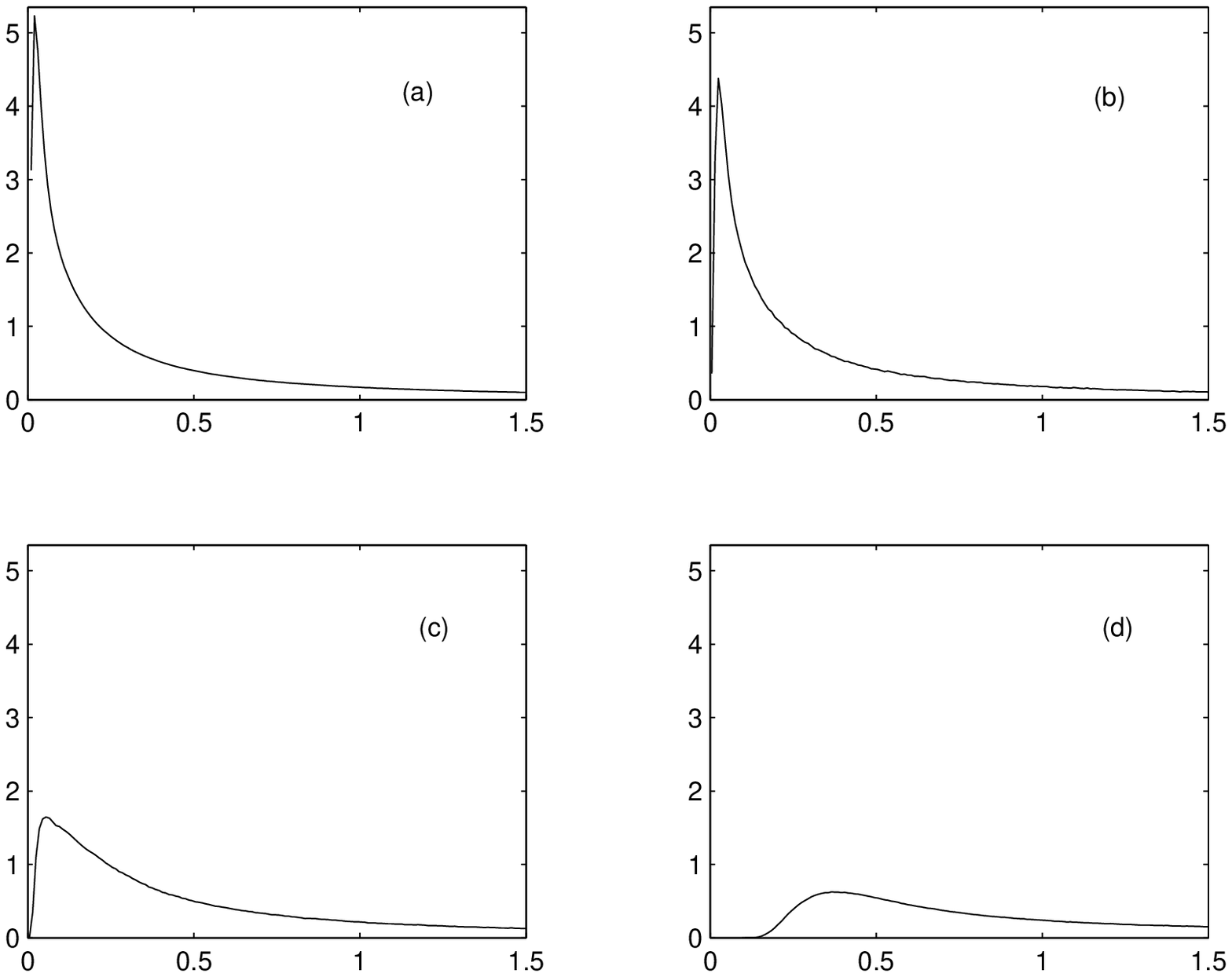}
			{10 cm}
			{ {\small Same as in Figura~\ref{fig:fig2} with $d=0.5$.}}		
	            {fig:fig3}
\end{figMacPc}
%
\section{Numerical comparisons\label{Section4}}
\hfill \\
The aim of this Section is to compare the behavior of the FPT
pdf's among Gauss-Markov processes and Gaussian non-Markov processes, in
order to analyze how the lack of memory affects the shape of the density,
also with reference to the specified type of correlation function.
For simplicity, set $x_0=0$ and $P[X(0)=0]=1,$ so that in the
following we shall consider the FPT pdf $g(t) := g[S(t),t|0,0]$.
\par
To be specific, we consider a stationary Gaussian process $X(t)$ with
zero mean and correlation function
\begin{equation}
\gamma(t):= e^{-\beta\, |t|}\, \cos (\alpha\, t), \qquad 
\beta \in { I\!\! R}^{+}
\label{(funcorr)}
\end{equation}
which is the simplest type of correlation having a concrete engineering
significance \cite{Yaglom87}. When the correlation function is of
type (\ref{(funcorr)}), $X(t)$
is  not mean square differentiable, since $\dot{\gamma}(0) \ne 0$.
Thus the series expansion (\ref{(0)}) does not hold. However, specific
assumptions on the parameter $\alpha$ help us characterize the shape of the FPT pdf.
\par
We start assuming $\alpha=0,$ so that the correlation function (\ref{(funcorr)}) can be
factorized as
$$
\gamma(t)= e^{-\beta\, \tau}\, e^{-\beta\,(t-\tau)} 
\qquad \beta \in {\bf R}^{+}, \, 0<\tau<t. 
$$
Hence, choosing $h_1(t)=e^{\beta t}$ and $h_2(t)=e^{-\beta t}$
in (\ref{eq:(2.2)}), $X(t)$ becomes
Gauss-Markov. Therefore, for any boundary $S(t),$ the FPT pdf $g(t)$ can be
numerically evaluated by solving the integral equation (\ref{eq:(3.1)}).
In the following, we consider boundaries of the form
\begin{equation}
S(t)=d\,e^{-\beta\, t}\, \Biggl\{ 1 - \frac{e^{2\, \beta\, t}-1}{2\, d^2}\, 
\ln \Biggl[ {1\over 4}+{1\over 4}\,
\sqrt{1+8 \exp\biggl(-{\ds\frac{4\,d^2}{e^{2\,\beta\, t} -1}\biggr)}}
\;\;\Biggr]\Biggr\},
\label{soglia}
\end{equation}
with $d>0$. It is evident that $\lim_{t\to 0}S(t)=d$ and that $S(t)$ tends to 0 as $t$ increases. 
Furthermore, as $d$ decreases, the boundary becomes flatter. In Figure~\ref{fig:fig1}, 
the boundary $S(t)$ given in (\ref{soglia}) is plotted
for $\beta=0.5$ and for two choices of the parameter $d$, i.e. $d=0.25$
and $d=0.5$. 
\par
As proved in \cite{DiNardo01bis} for boundaries of the form (\ref{soglia}), 
the FPT pdf $g(t)$ of a Gauss-Markov process admits the following closed form:
\begin{equation}
g(t) = {4 \, d \, \beta\, e^{\beta\, t}\over e^{2\,\beta\, t}-1}\;
\frac{\sqrt{1+8\, \exp\left(-{\ds\frac{4 d^2}{e^{2\,\beta\,t}
-1}}\right)}}{1+\sqrt{1+8\, \exp\left(-{\ds\frac{4 d^2}{e^{2\,\beta\, t}
-1}}\right)}}\;f[S(t),t|0,0],
\label{(FPTclosed)}
\end{equation}
where  $f[S(t),t|0,0]$ is the transition pdf of the Gauss-Markov process $X(t).$
\par
For a zero-mean Gauss-Markov process characterized by the  correlation function 
(\ref{(funcorr)}) with $\beta=0.5$ and $\alpha=0$, the FPT pdf $g(t)$ given by 
(\ref{(FPTclosed)}) through 
the boundary (\ref{soglia}) is plotted in Figure~2(a) for 
$d=0.25$ and in Figure~3(a) for $d=0.5$. Note that as $d$ increases  
the mode increase,  whereas the corresponding ordinate decrease.
\par
Setting $\alpha \neq 0$ in (\ref{(funcorr)}), the Gaussian process $X(t)$
is no longer Markov and its spectral density is given by
\begin{equation}
\Gamma(\omega) = \frac{2\, \beta\, (\omega^2+\alpha^2+\beta^2)}{\omega^4
+2\,\omega^2\,(\beta^2-\alpha^2)+(\beta^2+\alpha^2)^2}\,,
\label{(spectral)}
\end{equation}
thus being of a rational type. Since in (\ref{(spectral)}) the degree of the numerator is less 
than the degree of the denominator, it is possible to apply the
simulation algorithm described in Section 3 in order to estimate
the FPT pdf $\tilde{g}(t)$ of the process. 
\par
The simulation
procedure has been implemented by a parallel FORTRAN 90 code
on a 128-processor IBM SP4 supercomputer, based on MPI
language for parallel processing. The number of simulated
sample paths has been set equal to $10^7$. The estimated FPT pdf
$\tilde{g}(t)$ through the boundary (\ref{soglia}) with $\beta=0.5$ 
and $d=0.25$ are plotted in Figures~2(b)$\div$2(d) for Gaussian 
processes with correlation function (\ref{(funcorr)})
having $\alpha=10^{-10}, 0.25, 0.5$, respectively. 
For the same processes, Figures~3(b)$\div$3(d) show the estimated 
FPT pdf $\tilde{g}(t)$ through the boundary (\ref{soglia}) with $\beta=0.5$ 
and $d=0.5$.  Note that as $\alpha$ increases, the shape of the FPT pdf $\tilde{g}(t)$ 
becomes flatter and the related mode increases. Furthermore, as 
Figures 2(a)-2(b) and Figures 3(a)-3(b) show, $\tilde{g}(t)$  
is very similar to $g(t)$ for small values of $\alpha$. 

%
\ack
Work performed within a joint cooperation agreement between
Japan Science and Technology Corporation (JST) and Universit\`a di Napoli
Federico II, under partial support by INdAM (G.N.C.S).
%
%
%

\par\noindent
\noindent
{\em E.\ Di Nardo:} 
{\sc Dipartimento di Matematica, Universit\`a della Basilicata, Contrada Macchia Romana, Potenza, Italy}\\ 
{\em A.G.\ Nobile:} 
{\sc Dipartimento di Matematica e Informatica, Universit\`a di Sa\-lerno, 
Via S.\ Allende, I-84081  Baronissi (SA), Italy}\\ 
{\em E.\ Pirozzi and L.M.\ Ricciardi:} 
{\sc Dipartimento di Matematica e Applicazioni, Universit\`a di Napoli Federico II, 
Via Cintia, Napoli I-80126, Italy}

\end{document}